\documentclass[11pt,twoside, leqno]{article}

\usepackage{amssymb}
\usepackage{amsmath}

\allowdisplaybreaks

\def\rr{{\mathbb R}}
\def\rd{{{\rr}^d}}
\def\zz{{\mathbb Z}}

\def\cm{{\mathcal M}}

\def\nn{{\mathbb N}}
\def\dm{{\mathcal D}}

\def\cy{{\mathcal Y}}

\def\fz{\infty}
\def\az{\alpha}
\def\supp{{\mathop\mathrm{\,supp\,}}}

\def\loc{{\mathop\mathrm{\,loc\,}}}

\def\lz{\lambda}
\def\dz{\delta}

\def\ez{\epsilon}

\def\bz{\beta}

\def\gz{{\gamma}}

\def\vz{\varphi}

\def\wz{\widetilde}

\def\hs{\hspace{0.3cm}}

\def\ls{\lesssim}

\def\paz{{\partial}}

\def\lhp{{h_{\rm atb}^{1,\,p}(\mu)}}

\def\lhg{{h_{{\rm atb},\,\gz}^{1,\,p}(\mu)}}

\def\hp{{H_{\rm atb}^{1,\,p}(\mu)}}
\def\hpg{{H_{{\rm atb},\,\gz}^{1,\,p}(\mu)}}

\def\lo{{L^1(\mu)}}
\def\wlo{{L^{1,\,\fz}(\mu)}}
\def\lt{{L^2(\mu)}}
\def\lin{{L^\fz(\mu)}}
\def\lp{{L^p(\mu)}}
\def\lho{{h^1(\mu)}}

\def\bmo{{\mathop\mathrm{RBMO}\,(\mu)}}
\def\lbmo{{\mathop\mathrm{rbmo}\,(\mu)}}

\def\dint{\displaystyle\int}

\def\dfrac{\displaystyle\frac}
\def\dsup{\displaystyle\sup}

\def\r{\right}
\def\lf{\left}

\newtheorem{thm}{Theorem}[section]

\newtheorem{prop}{Proposition}[section]
\newtheorem{rem}{Remark}[section]

\newtheorem{defn}{Definition}[section]

\newtheorem{pf}{\bf Proof.}

\numberwithin{equation}{section}

\begin{document}

\arraycolsep=1pt

\title{{\vspace{-5cm}\small\hfill\bf Georgian Math. J., to appear}\\
\vspace{4cm}\Large\bf Boundedness of Linear Operators via Atoms on Hardy Spaces
with Non-doubling Measures\footnotetext{\hspace{-0.35cm} 2000 {\it Mathematics Subject
Classification}. {Primary 42B20; Secondary 42B30, 42B35.}
\endgraf{\it Key words and phrases.} non-doubling measure,
$H^1(\mu)$, $\lho$, linear operator,
atomic block, block, Calder\'on-Zygmund operator,
fractional integral operator, commutator.
\endgraf
The first author is supported by the National
Natural Science Foundation (Grant No. 10871025) of China.}}
\author{Dachun Yang and Dongyong Yang}
\date{ }
\maketitle

\begin{center}
\begin{minipage}{13.5cm}\small
{\noindent{\bf Abstract.} Let $\mu$ be a non-negative Radon
measure on ${\mathbb R}^d$ which only satisfies the polynomial growth
condition. Let ${\mathcal Y}$ be a Banach space
and $H^1(\mu)$ the Hardy space of Tolsa. In this paper,
the authors prove that a linear operator $T$ is
bounded from $H^1(\mu)$ to ${\mathcal Y}$
if and only if $T$ maps all $(p, \gamma)$-atomic blocks into uniformly bounded
elements of ${\mathcal Y}$; moreover, the authors prove that for
a sublinear operator $T$ bounded from $L^1(\mu)$ to $L^{1,\,\infty}(\mu)$,
if $T$ maps all $(p, \gamma)$-atomic blocks
with $p\in(1, \infty)$ and $\gamma\in{\mathbb N}$
into uniformly bounded elements of $L^1(\mu)$,
then $T$ extends to a bounded sublinear operator from $H^1(\mu)$
to $L^1(\mu)$. For the localized atomic
Hardy space $h^1(\mu)$, corresponding results are also
presented. Finally, these results are applied to Calder\'on-Zygmund operators,
Riesz potentials and multilinear commutators generated by Calder\'on-Zygmund operators
or fractional integral operators with Lipschitz functions, to
simplify the existing proofs in the corresponding papers.}
\end{minipage}
\end{center}

\section{Introduction\label{s1}}

\hskip\parindent The real-variable theory of Hardy spaces on $\rd$, which began
with the remarkable work of Stein and Weiss \cite{sw60}, has been
transformed into a rich theory. The well-known atomic
and molecular characterizations of Hardy spaces
enable one to deduce the boundedness on Hardy spaces of (sub)linear operators
from their behaviors on atoms or molecules in principle. However, Meyer,
Taibleson and Weiss \cite{mtw85}
constructed an example of $f\in H^1(\rd)$ such that its
norm can not be achieved by its finite atomic decompositions via
$(1, \fz)$-atoms. Inspired by this, Bownik \cite{b05} showed
that there exists a linear functional,
which maps all $(1,\fz)$-atoms of $H^1(\rd)$ into bounded scalars but
does not admit a bounded extension to $H^1(\rd)$. It turns out that
the condition that a linear operator $T$ maps all $(1, \fz)$-atoms
into a uniformly bounded subset of certain
quasi-Banach space ${\mathcal B}$ fails to guarantee the extension
of $T$ to a bounded linear operator from the whole $H^1(\rd)$ to ${\mathcal B}$.
Recently, Meda, Sj\"ogren and Vallarino \cite{msv08} proved that any
linear operator mapping all $(1,q)$-atoms with $q\in(1, \fz)$ or all
continuous $(1,\fz)$-atoms into a uniformly bounded elements in a
given Banach space ${\mathcal B}$ extends to a bounded linear operator from
$H^1(\rd)$ to ${\mathcal B}$. Independently, in \cite{yz08}, a boundedness
criterion was established as follows: a non-negative sublinear operator $T$ extends
to a bounded sublinear operator from Hardy spaces $H^p(\rd)$ with $p\in
(0, 1]$ to certain quasi-Banach space ${\mathcal B}$ if and only if $T$
maps all $(p, 2)$-atoms into uniformly bounded elements of
${\mathcal B}$.
On the other hand, via making clear
the dual and the completion of the space of finite
linear combinations of $(p, \fz)$-atoms with $p\in(0, 1]$,
Ricci and Verdera \cite{rv} further proved that if $T$ is a linear operator
mapping all $(p, \fz)$-atoms with $p\in(0, 1)$ uniformly
bounded to a Banach space ${\mathcal B}$, then $T$ extends
to a bounded linear operator from $H^p(\rd)$ to ${\mathcal B}$.

Let $\mu$ be a non-negative Radon
measure on ${\mathbb R}^d$ which only satisfies the polynomial growth
condition. Let ${\mathcal Y}$ be a Banach space
and $H^1(\mu)$ the Hardy space of Tolsa (see \cite{t01,t03}). In this paper,
we prove that a linear operator $T$ is
bounded from $H^1(\mu)$ to ${\mathcal Y}$
if and only if $T$ maps all $(p, \gamma)$-atomic blocks
(\cite{t01, hmy05}) into uniformly bounded
elements of ${\mathcal Y}$; moreover, we show that for
a sublinear operator $T$ bounded from $L^1(\mu)$ to $L^{1,\,\infty}(\mu)$,
if $T$ maps all $(p, \gamma)$-atomic blocks
with $p\in(1, \infty)$ and $\gamma\in{\mathbb N}$
into uniformly bounded elements of $L^1(\mu)$,
then $T$ extends to a bounded sublinear operator from $H^1(\mu)$
to $L^1(\mu)$. For the localized atomic
Hardy space $h^1(\mu)$ in \cite{hyy07}, corresponding results are also
presented. Finally, these results are applied to Calder\'on-Zygmund operators,
Riesz potentials and multilinear commutators generated by Calder\'on-Zygmund operators
or fractional integral operators with Lipschitz functions, to
simplify the existing proofs in the corresponding papers \cite{cmy05,hmy08,my06}.
Moreover, these results seal a gap existing in the proof
of \cite[Theorem 1.1]{hmy08}.

Recall that a non-negative Radon measure $\mu$ on $\rd$ is called
a non-doubling measure, if there exist positive constants
$C$ and $n\in (0,d]$ such that for all $x\in\rd$ and $r>0$,
\begin{equation*}
\mu\lf(B(x,r)\r)\le Cr^n,
\end{equation*}
where $B(x,r)\equiv \{y\in\rd: |x-y|<r\}.$
 Such a measure $\mu$ is not necessary to be doubling, which is
a crucial assumption in the classical theory of harmonic analysis. In
recent years, it was shown that many classical results concerning the
theory of Calder\'on-Zygmund operators and function spaces
remain valid for non-doubling measures; see, for example,
\cite{ntv97,ntv98,t01,t01b,t01c,ntv02,ntv03}.
Moreover, the harmonic analysis for non-doubling measures plays
an important role in the solution of several long-standing open questions
related to analytic capacity, like Painlev\'e's problem
and Vitushkin's conjecture; see \cite{t03b,t05,v02,v03}
for more details.

To state the main results of this paper, we first recall some notation and
notions.

Throughout this paper, by a cube $Q\subset\rd$,
we mean a closed cube whose sides are parallel to the axes and
centered at certain point of $\supp(\mu)$, and we denote its side
length by $l(Q)$ and its center by $x_Q$.
For any given $\lambda\in(0, \fz)$ and cube $Q$, $\lambda Q$ denotes the
cube concentric with $Q$ and having side length $\lz l(Q)$.
Given two cubes $Q,\ R\subset\rd$,
let $Q_R$ be the smallest cube concentric with $Q$ containing $Q$
and $R$. We also set
$\nn\equiv\{1,\,2,\,\cdots\}$ and $\zz_+\equiv\nn\cup\{0\}$.

The following coefficient was first introduced
by Tolsa in \cite{t01} and the Hardy space $H^1(\mu)$
by Tolsa in \cite{t03}.

\begin{defn}\label{d1.1}\rm
Given two cubes $Q,\, R \subset \rd$, define
$$\dz(Q,R)\equiv\max\lf\{\dint_{Q_R\setminus Q }\frac1{|x - x_Q|^n }\,d\mu(x),\
\dint_{R_Q\setminus R} \frac1{|x- x_R|^n}\, d\mu(x)\r\}.$$
\end{defn}

\begin{defn}\label{d1.2}\rm
 Given $f\in L^1_{\loc}(\mu)$, set
$$\cm_\Phi (f)(x)\equiv\dsup_{\vz\sim x}\lf|\dint_\rd f\vz\,d\mu\r|,$$
where the notation $\varphi\sim x$ means that $\varphi\in\lo\cap
C^1(\rd)$ and  satisfies
\begin{itemize}
    \item [{\rm (i)}] $\|\varphi\|_\lo\le 1,$
    \item [{\rm (ii)}]$0\le \varphi(y)\le\frac1{|y-x|^n}$ for all $y\in\rd$, and
    \item [{\rm (iii)}]$|\nabla\varphi(y)|\le\frac1{|y-x|^{n+1}}$ for all $y\in \rd$,
  where $\nabla=(\frac {\paz}{\paz x_1},\cdots,\frac {\paz}{\paz x_d})$.
\end{itemize}
\end{defn}

\begin{defn}\label{d1.3}\rm
The Hardy space $H^1(\mu)$ is defined to be the set of all functions
$f\in \lo$ satisfying that $\int_\rd f\,d\mu=0$ and $\cm_\Phi
(f)\in\lo$. Moreover, the norm of $f\in H^1(\mu)$ is defined by
$$\|f\|_{H^1(\mu)}\equiv\|f\|_\lo+\|\cm_\Phi (f)\|_\lo.$$
\end{defn}

We now recall atomic characterizations of the Hardy space $H^1(\mu)$
and its localized variant in \cite{t01, t03, hmy05, hyy07}.

\begin{defn}\label{d1.4}\rm
Let $\eta \in(1,\fz)$, $\gz\in\nn$ and $p\in(1,\fz]$.
A function $b \in L^1_{\loc}(\mu)$ is called
a $(p, \gz)$-atomic block if
\begin{enumerate}
\item [{\rm(i)}]there exists certain cube $R$ such that $\supp (b) \subset R$,
\item [{\rm(ii)}]$\int_\rd b(x) \,d\mu(x) = 0$,
\item [{\rm(iii)}]for $j = 1, 2$, there exist functions $a_j$
supported on cubes $Q_j\subset R$ and numbers $\lz_j\in\rr$ such
that $b = \lz_1a_1 + \lz_2a_2$, and
$$\|a_j\|_\lp\le  [\mu(\eta Q_j )]^{1/p-1} [1+\dz(Q_j,\,R)]^{-\gz}.$$
\end{enumerate}
Then we define
$|b|_{H^{1,\,p}_{{\rm atb},\,\gz}(\mu)}\equiv|\lz_1| + |\lz_2|.$

A function $f\in\lo$ is said to belong to the space  $H^{1,\,p}_{{\rm atb},\,\gz}(\mu)$
if there exist $(p, \gz)$-atomic blocks $\{b_i\}_{i\in\nn}$ such that
$f =\sum_{i=1}^\fz b_i$
with $\sum^\fz_{i=1} |b_i|_\hpg<\fz$. The $\hpg$ norm of $f$
is defined by
$$\|f\|_\hpg\equiv \inf\lf\{\sum_{i=1}^\fz|b_i|_\hpg\r\},$$
where the infimum is taken over all the possible decompositions of
$f$ as above.
\end{defn}

\begin{rem}\label{r1.1}\rm If $\gz=1$, we denote $\hpg$ simply by $\hp$.
The space $\hpg$ when $\gz=1$ was introduced by Tolsa in \cite{t01}, and
when $\gz>1$ was introduced in \cite{hmy05}.
It was proved in \cite{t01, t03, hmy05} that the definition of $\hpg$
is independent of the chosen constant $\eta\in(1,\fz)$
and that all the atomic Hardy spaces $\hpg$ with
$\gz\in\nn$ and $p\in(1,\fz]$ coincide with
$H^1(\mu)$ with equivalent norms.
In the rest of this paper, unless explicitly stated,
we always choose $\eta=2$  and $\gz=1$ in the definition of $\hpg$.
\end{rem}

We now recall the notions of initial cubes and
the localized atomic Hardy space, respectively, in \cite{t01b} and
\cite{hyy07}.

\begin{defn}\label{d1.5}\rm
The Euclidean space $\rd$ is called an initial cube
if $\dz(Q,\rd) <\fz$ for certain cube $Q$ with $l(Q)\in (0, \fz)$.
\end{defn}

\begin{rem}\label{r1.2}\rm
In \cite[p.\,67]{t01b}, it was pointed out that if $\dz(Q, \rd)<\fz$
for certain cube $Q$ with $l(Q)\in (0,\fz)$, then $\dz(Q', \rd)<\fz$ for any
cube $Q'$ with $l(Q')\in (0,\fz)$.
\end{rem}

Let $A$ be a big positive constant. In particular, as in
\cite{t01b, t03}, we
assume that $A$ is much bigger than the constant $\ez_1$
in Lemma 3.2 of \cite{t01b}.
In the case that $\rd$ is not an initial cube, let
$\{R_{-j}\}_{j\in\zz_{+}}$ be a sequence of increasing concentric
`reference' cubes as in \cite{t01b} and
\begin{eqnarray*}
\dm&\equiv&\big\{Q\subset\rd:\ \mathrm{there\ exists\ a\ cube}\
P\subset Q\ \mathrm{and}\ j\in \zz_{+}\
\mathrm{such\ that}\\
&&\hs P\subset R_{-j}\ \mathrm{with}\ \dz(P, R_{-j})\le (j+1)A+\ez_1\big\}.
\end{eqnarray*}
If $\rd$ is an initial cube, we then define the set
$$\dm\equiv\big\{Q\subset\rd:\ \mathrm{there\
exists\ a\ cube}\ P\subset Q\ \mathrm{such\ that}\ \dz(P, \rd)\le A+\ez_1\big\}.$$
It  was pointed out in \cite{hyy07} that the definition of the set $\dm$
is independent of the chosen reference cubes $\{R_{-j}\}_{j\in\zz_{+}}$ in the
sense modulo certain small error; see also \cite[p.\,68]{t01b}.

\begin{defn}\label{d1.6}\rm
Let $\eta\in(1,\fz)$, $\gz\in\nn$ and $p\in(1, \fz]$.
A function $b\in L^1_{\loc}(\mu)$ is called
a $(p, \gz)$-block if only (i) and
(iii) of Definition \ref{d1.4}  hold.
Moreover, define $|b|_{\lhg}\equiv\sum_{j=1}^2|\lz_j|$.

A function $f\in\lo$ is said to belong to the space $\lhg$
if there exist $(p, \gz)$-atomic blocks or $(p, \gz)$-blocks $\{b_i\}_i$ such that
$f=\sum_i b_i$ and $\sum_i|b_i|_{\lhg}<\fz$, where $b_i$ is a
$(p, \gz)$-atomic block as in Definition \ref{d1.4}
if $\supp(b_i)\subset R_i$ and $R_i\notin\dm$, while
$b_i$ is a $(p, \gz)$-block if $\supp(b_i)\subset R_i$ and $R_i\in\dm$.
Moreover, the $\lhg$ norm of $f$ is defined by
$$\|f\|_\lhg\equiv\inf\lf\{\sum_i|b_i|_\lhg\r\},$$
where the infimum is taken over all decompositions of $f$ as above.
\end{defn}

\begin{rem}\label{r1.3}\rm
When $\gz=1$, we denote the space $\lhg$ simply by $\lhp$,
which was introduced in \cite{hyy07};
moreover, it was proved there that the definition of
$h^{1,\,p}_{\rm atb}(\mu)$
is independent of the chosen constant $\eta\in(1,\fz)$, and that
all the localized atomic Hardy spaces
$h^{1,\,p}_{\rm atb}(\mu)$ with $p\in(1,\fz)$ coincide with
$h^{1,\,\fz}_{\rm atb}(\mu)$ with equivalent norms.
\end{rem}

By the same argument as that used in the proof of Theorem 2.1
in \cite{hmy05}, we have the following equivalent atomic characterization of
$\lhg$. We omit the details here.

\begin{prop}\label{p1.1}
Let $\eta\in(1, \fz)$, $\gz\in\nn$ with $\gz>1$ and $p\in (1, \fz]$. Then
$\lhg=h^{1,\,p}_{\rm atb}(\mu)$ with equivalent norms.
\end{prop}

As a consequence of Remark \ref{r1.3} and Proposition \ref{p1.1},
throughout this paper,  we
denote $\lhg$ simply by $h^1(\mu)$. Moreover, unless
explicitly stated, in what follows,
we always choose $\eta=2$ and $\gz=1$ in the definition of $\lhg$.

The main results of this paper are as follows.

\begin{thm}\label{t1.1}
Let $\eta\in(1,\fz)$, $\gz\in\nn$, $p\in(1, \fz)$, $T$ be a linear operator
and ${\mathcal Y}$ a Banach space.

(i) If there exists a non-negative constant $C$ such that for all
$(p, \gz)$-atomic blocks $b$,
\begin{equation}\label{1.1}
||Tb\|_{\mathcal Y}\le C |b|_{\hpg},
\end{equation}
then $T$ extends to a bounded linear operator
from $H^1(\mu)$ to $\cy$.

(ii) If there exists a non-negative constant $\wz C$ such that for all
$(p, \gz)$-atomic blocks $b$ with $\supp(b)\subset R$ and
$R\notin\dm$, and all $(p, \gz)$-blocks $b$ with $\supp(b)\subset R$ and $R\in\dm$,
\begin{equation}\label{1.2}
||Tb\|_{\mathcal Y}\le \wz C |b|_{\lhg},
\end{equation}
then $T$ extends to a bounded linear operator from $\lho$ to $\cy$.
\end{thm}

\begin{rem}\label{r1.4}\rm
Observe that \eqref{1.1} (or \eqref{1.2})
is also necessary for an operator $T$ to be
bounded from $H^1(\mu)$ (or $\lho$) to $\cy$.
From this fact and Theorem \ref{t1.1}, we further deduce that
if $T$ is linear, then $T$ extends to a bounded linear operator
from $H^1(\mu)$ (or $\lho$) to $\cy$
if and only if $T$ satisfies \eqref{1.1} (or \eqref{1.2}).
\end{rem}

For sublinear operators bounded from $\lo$ to $\wlo$,
we also have the following conclusion.

\begin{thm}\label{t1.2}
Let $\eta\in(1,\fz)$, $\gz\in\nn$,
$p\in(1, \fz)$ and $T$ be a sublinear operator bounded from $\lo$ to $\wlo$.

(i)  If $T$ satisfies \eqref{1.1} with $\cy=L^1(\mu)$, then $T$ extends to a bounded
sublinear operator from $H^1(\mu)$ to $\lo$.

(ii) If $T$ satisfies \eqref{1.2} with $\cy=L^1(\mu)$, then $T$
extends to a bounded sublinear operator
from $\lho$ to $\lo$.
\end{thm}

Proofs of Theorems \ref{t1.1} and \ref{t1.2} are given in
Section \ref{s2}. We remark that the proof of Theorem \ref{t1.2}
would be trivial if $T$ were linear. In fact, it is
easy to see that if the linear operator $T$ is
continuous from $L^1(\mu)$ to $L^{1,\,\fz}(\mu)$,
and the image of atomic blocks (or blocks) is uniformly
bounded in $L^1(\mu)$, then $T$ is automatically
bounded from $H^1(\mu)$ (or $h^1(\mu)$) to $L^1(\mu)$.
For sublinear operators, the proof of Theorem \ref{t1.2}
requires only an easy additional measure theoretic
argument.

In Section \ref{s3},  we apply Theorem \ref{t1.1} to
Calder\'on-Zygmund operators, Riesz potentials and multilinear commutators
generated by Calder\'on-Zygmund operators
or fractional integral operators with Lipschitz functions,
to simplify the existing proofs in the corresponding papers;
see \cite[Theorem 1]{cmy05}, \cite[Theorem 1.1]{hmy08}
and \cite[Theorems 3.1, 4.2]{my06}. In particular,
we seal a gap existing in the proof that (III) implies (IV)
of \cite[Theorem 1.1]{hmy08} (see \cite[pp.\,379-381]{hmy08}).
We also prove that if $\rd$ is an initial cube,
then the Calder\'on-Zygmund operator is bounded from $\lho$ to $\lo$.

We now make some conventions. Throughout this paper, we always use
$C$ to denote a positive constant that is independent of the main
parameters involved but whose value may differ from line to line.
Constants with subscripts, such as $C_1$, do not change in
different occurrences. If $f\le Cg$, we then write $f\ls g$;
and if $f \ls g\ls f$, we write $f\sim g.$

\section{Proofs of Theorems \ref{t1.1} and \ref{t1.2}\label{s2}}

\hskip\parindent
In this section, we show Theorems \ref{t1.1} and \ref{t1.2}.
To start with, we recall some useful notions and notation.

Let $p\in(1, \fz]$, $L^p_c(\mu)$ be the space of functions in
$\lp$ with compact support and $L^p_{c,\,0}(\mu)$
the space of functions in
$L^p_c(\mu)$ having integral 0.
Moreover, for each cube $Q$, we denote by $L^p(Q)$ the subspace
of functions in $L^p(\mu)$ supported in $Q$ and $L^p_0(Q)
\equiv L^p_{c,\,0}(\mu)\cap L^p(Q)$.
Then the unions of $L^p_0(Q)$ and $L^p(Q)$ as $Q$ varies over all cubes
coincide with $L^p_{c,\,0}(\mu)$ and $L^p_c(\mu)$, respectively.
Now let $\{Q_j\}_{j\in\nn}$ be a sequence of increasing concentric cubes
with $\rd=\cup_{j\in\nn} Q_j$. We topologize $L^p_{c,\,0}(\mu)$ (resp. $L^p_c(\mu)$) as
the strict inductive limit of the spaces $L^p_0(Q_j)$ (resp. $L^p(Q_j)$)
(see \cite[II,\, p.\,33]{b87} for the definition of the
strict inductive limit topology).
It is known that the definition of the topology of
$L^p_{c,\,0}(\mu)$ (resp. $L^p_c(\mu)$) is independent of the choice of
$\{Q_j\}_{j\in\nn}$.

We now recall the definitions of $\bmo$ of Tolsa in \cite{t01}
and $\lbmo$ in \cite{hyy07}.

\begin{defn}\label{d2.1}\rm
(i) Let $p\in[1, \fz)$.
A function $f\in L^1_\loc(\mu)$ is said to be in the
space $\bmo$ if there exists a nonnegative constant $C$ such
that for any doubling cube $Q$,
\begin{equation}\label{2.1}
\lf[\dfrac 1{\mu(Q)}\dint _{Q}\lf|f(y)-m_Q(f)\r|^p\,d\mu(y)\r]^{1/p}\le C,
\end{equation}
and for any two doubling cubes $Q\subset R$,
\begin{equation}\label{2.2}
|m_Q(f)-m_R(f)|\le C[1+\dz(Q, R)],
\end{equation}
where $m_Q(f)$ denotes the mean of $f$ over cube $Q$,
namely, $m_Q(f)\equiv\frac 1{\mu(Q)}\int_Qf(y)\,d\mu(y)$. Moreover, we
define the $\bmo$ norm of $f$ to be the minimal constant
$C$ as above and denote it by $\|f\|_\bmo$.

(ii) Let $p\in[1, \fz)$.
A function $f\in L^1_\loc(\mu)$ is said to be in the
space $\lbmo$ if there exists a nonnegative constant $C$ such
that \eqref{2.1} holds for any doubling cube $Q\notin \dm$,
\eqref{2.2} holds for any two doubling cubes $Q\subset R$ with $Q\notin \dm$,
and  for any doubling cube $Q\in\dm$,
\begin{equation*}
\lf[\dfrac 1{\mu(Q)}\dint _{Q}|f(y)|^p\,d\mu(y)\r]^{1/p}\le C.
\end{equation*}
Moreover, we define the $\lbmo$ norm of $f$ to be the minimal constant
$C$ as above and denote it by $\|f\|_\lbmo$.
\end{defn}

\begin{rem}\label{r2.1}\rm
In \cite{t01}, Tolsa showed that $\bmo$ is the dual space of $H^1(\mu)$.
On the other hand, it was proved in \cite{hyy07}
that $\lbmo$ is the dual space of $\lho$.
\end{rem}

\newtheorem{pftf}{\bf  Proof of Theorem \ref{t1.1}.}
\renewcommand\thepftf{}

\begin{pftf}\rm
We first show (i) of Theorem \ref{t1.1}. To this end, without loss of generality,
we may assume $p=2$. Moreover, by Remark \ref{r1.1},
we choose $\eta=2$ and $\gz=1$ in the definition of $\hpg$.
Let $Q$ be a fixed cube.
If $f\in L^2_0(Q)$, then $f$ is a $(2, 1)$-atomic block and
\begin{equation}\label{2.3}
|f|_{H^{1,\,2}_{\rm \,atb}(\mu)}\le \|f\|_\lt[\mu(2Q)]^{1/2}.
\end{equation}
Moreover, from this and \eqref{1.1}, it follows that for any
sequence of increasing concentric cubes $\{Q_j\}_{j\in\nn}$
with $\rd=\cup_{j\in\nn} Q_j$,
$T$ is bounded from $L^2_0(Q_j)$ to $\cy$ for each $j\in\nn$.
Then $T$ is bounded from $L^2_{c,\,0}(\mu)$ to $\cy$, which implies that
the adjoint operator $T^\ast$ of $T$ is bounded from the dual space
$\cy^\ast$ of $\cy$ to $[L^2_{c,\,0}(\mu)]^\ast$.
Moreover, for all functions $f\in\cy^\ast$ and $(2, 1)$-atomic blocks $b$,
we have
\begin{equation}\label{2.4}
\lf|\dint_\rd b(x)T^\ast(f)(x)\,d\mu(x)\r|
=|\langle Tb, f\rangle|\ls
\|f\|_{\cy^\ast}|b|_{H^{1,\,2}_{\rm atb}(\mu)}.
\end{equation}

We claim that for all $f\in\cy^\ast$,
$T^\ast f\in\bmo$ and $\|T^\ast f\|_\bmo\ls \|f\|_{\cy^\ast}$.
In fact, observe that for any doubling cube $Q$ and
$\phi\in L^2(Q)$ with $\|\phi\|_{L^2(Q)}=1$, $[\phi-m_Q(\phi)]\chi_Q$ is
a $(2, 1)$-atomic block, where and in what follows,
$\chi_Q$ denotes the characteristic function of the set
$Q$. From this, \eqref{2.3} and \eqref{2.4}, we deduce that
\begin{eqnarray*}
\lf[\dint_Q\lf|T^\ast f(x)-m_Q(T^\ast f)\r|^2\,d\mu(x)\r]^{1/2}
&=&\dsup_{\|\phi\|_{L^2(Q)}=1}\lf|\dint_Q\phi(x)
\lf[T^\ast f(x)-m_Q(T^\ast f)\r]\,d\mu(x)\r|\\
&=&\dsup_{\|\phi\|_{L^2(Q)}=1}\lf|\dint_Q[\phi(x)-m_Q(\phi)]T^\ast f(x)\,d\mu(x)\r|\\
&\ls&\|f\|_{\cy^\ast}[\mu(Q)]^{1/2},
\end{eqnarray*}
which implies that
\begin{equation}\label{2.5}
\lf[\dfrac1{\mu(Q)}\dint_Q\lf|T^\ast f(x)-m_Q(T^\ast f)\r|^2
\,d\mu(x)\r]^{1/2}\ls \|f\|_{\cy^\ast}.
\end{equation}
By \eqref{2.5} and Definition \ref{d2.1} (i), the claim is reduced to
showing that for all doubling cubes $Q\subset R$,
\begin{equation}\label{2.6}
\lf|m_Q(T^\ast f)-m_R(T^\ast f)\r|\ls [1+\dz(Q, R)]\|f\|_{\cy^\ast}.
\end{equation}
Let
\begin{equation*}
a_1\equiv\frac{|T^\ast f-m_{R}(T^\ast f)|^{2}}
{T^\ast f-m_{R}(T^\ast f)}\chi_{Q\cap\{T^\ast f\not=m_{R}(T^\ast f)\}},
\end{equation*}
$a_2\equiv C_R\chi_R$ and $b\equiv a_1+a_2$, where $C_R$ is a constant such that
$b$ has integral 0. Then $b$ is a $(2, 1)$-atomic block and
\begin{eqnarray*}
|b|_{H^{1,\,2}_{\rm atb}(\mu)}&\ls& \|a_1\|_{L^2(\mu)}[\mu(Q)]^{1/2}[1+\dz(Q, R)]
+|C_R|\mu(R)\\
&\ls&\lf[\dint_Q|T^\ast f(x)-m_R(T^\ast f)|^{2}\,d\mu(x)\r]^{1/2}
[\mu(Q)]^{1/2}[1+\dz(Q, R)].
\end{eqnarray*}
By this, \eqref{2.4} and \eqref{2.5} with $Q$ replaced by $R$, we have
\begin{eqnarray*}
&&\dint_Q|T^\ast f(x)-m_{R}(T^\ast f)|^{2}\,d\mu(x)\\
&&\hs=\dint_\rd a_1(x)[T^\ast f(x)-m_{R}(T^\ast f)]\,d\mu(x)\\
 && \hs\le\lf[\Big|\dint_\rd b(x)T^\ast f(x)\,d\mu(x)\Big|
 +|C_R|\dint_R|T^\ast f(x)-m_{R}(T^\ast f)|\,d\mu(x)\r]\\
 &&\hs\ls\lf[\dint_Q|T^\ast f(x)-m_{R}(T^\ast f)|^{2}\,d\mu(x)\r]^{1/2}
[\mu(Q)]^{1/2}[1+\dz(Q, R)]\|f\|_{\cy^\ast},
\end{eqnarray*}
which implies that
$$\lf[\frac1{\mu(Q)}\dint_Q|T^\ast f(x)-m_{R}(T^\ast f)|^{2}\,d\mu(x)\r]^{1/2}
\ls[1+\dz(Q, R)]\|f\|_{\cy^\ast}.$$
From this, the H\"older inequality and
\eqref{2.5}, it then follows that
\begin{eqnarray*}
|m_Q(T^\ast f)-m_R(T^\ast f)|&\le&\frac1{\mu(Q)}
\dint_Q\lf[|m_Q(T^\ast f)-T^\ast f(x)|+
|T^\ast f(x)-m_R(T^\ast f)|\r]\,d\mu(x)\\
&\ls& [1+\dz(Q, R)]\|f\|_{\cy^\ast},
\end{eqnarray*}
which implies \eqref{2.6}. By this together with \eqref{2.4},
we obtain that $T^\ast f\in\bmo$ and
$\|T^\ast f\|_\bmo\ls \|f\|_{\cy^\ast}$. Thus, the claim is true.

Let $H^{1,\,2}_{\rm fin}(\mu)$ be the set of all finite
linear combinations of $(2, 1)$-atomic blocks. Then $H^{1,\,2}_{\rm fin}(\mu)$
is dense in $H^1(\mu)$. On the other hand, $H^{1,\,2}_{\rm fin}(\mu)$
coincides with $L^2_{c,\,0}(\mu)$ as vector spaces. Then
by Remark \ref{r2.1} and the above claim, we
have that for all $g\in H^{1,\,2}_{\rm fin}(\mu)$ and $f\in
\cy^\ast$ with $\|f\|_{\cy^\ast}=1$,
$\lf|\langle Tg, f\rangle\r|=\lf|\langle g, T^\ast f\rangle\r|\ls
\|g\|_{H^1(\mu)}\|T^\ast f\|_\bmo\ls \|g\|_{H^1(\mu)}.$
From this  and \eqref{1.1}, it follows that $Tg\in\cy$ and
$\|Tg\|_\cy\ls \|g\|_{H^1(\mu)},$
which via a density argument then completes
the proof of Theorem \ref{t1.1} (i).

We now prove (ii). Similarly to (i), without loss of generality,
we may assume that $p=2$ and we choose $\eta=2$ and $\gz=1$
in the definition of $\lhg$.
Using an argument similar to (i), we see that if $T$
satisfies \eqref{1.2}, then
$T$ is bounded from $L^2_c(\mu)$ to $\cy$, which implies that
$T^\ast$ is bounded from $\cy^\ast$ to $[L^2_c(\mu)]^\ast$.
Moreover, we have that for all $f\in\cy^\ast$ and $(2, 1)$-atomic blocks
or $(2, 1)$-blocks $b$ as in \eqref{1.2},
\begin{equation}\label{2.7}
\lf|\dint_\rd b(x)T^\ast(f)(x)\,d\mu(x)\r|
=|\langle Tb, f\rangle|\ls
\|f\|_{\cy^\ast}|b|_{h^{1,\,2}_{\rm atb}(\mu)}.
\end{equation}

We claim that for all $f\in\cy^\ast$, $T^\ast f\in \lbmo$
and $\|T^\ast f\|_\lbmo\ls\|f\|_{\cy^\ast}$. In fact,
we first prove that for any doubling cube $Q\in\dm$,
\begin{equation}\label{2.8}
|m_Q(T^\ast f)|\ls \|f\|_{\cy^\ast}.
\end{equation}
Let $Q\in\dm$ be doubling. Observe that
for any doubling cube $Q$ and
$\phi\in L^2(Q)$ with $\|\phi\|_{L^2(Q)}=1$, $\phi$ is
a $(2, 1)$-block. From this and \eqref{2.7}, it follows  that
\begin{eqnarray*}
\lf[\dint_Q\lf|T^\ast f(x)\r|^2\,d\mu(x)\r]^{1/2}
&=&\dsup_{\|\phi\|_{L^2(Q)}=1}\lf|\dint_Q\phi(x)T^\ast f(x)\,d\mu(x)\r|\\
&\ls&\|f\|_{\cy^\ast}|\phi|_{h^{1,\,2}_{\rm atb}(\mu)}
\ls\|f\|_{\cy^\ast}[\mu(Q)]^{1/2},
\end{eqnarray*}
which via  the H\"older inequality yields \eqref{2.8}.

By the proof of \eqref{2.5}, we also have that for any doubling cube $Q\notin \dm$,
\begin{equation}\label{2.9}
\frac1{\mu(Q)}\dint_Q\lf|T^\ast f(x)-m_Q(T^\ast f)\r|\,d\mu(x)
 \ls \|f\|_{\cy^\ast}.
 \end{equation}
By this and \eqref{2.8} together with Definition \ref{d2.1} (ii),
to show the claim, it suffices to prove
that for any two doubling cubes $Q\subset R$ with
$Q\notin\dm$,
\begin{equation}\label{2.10}
|m_Q(T^\ast f)-m_R(T^\ast f)|\ls [1+\dz(Q, R)]\|f\|_{\cy^\ast}.
\end{equation}
In fact, if $R\notin\dm$, then by the proof
of \eqref{2.6}, we obtain \eqref{2.10}. Now suppose that
$R\in \dm$. We set $a\equiv a_1$, where $a_1$ is as in the proof of \eqref{2.6}.
Then $a$ is a $(2, 1)$-block with $\supp(a)\subset R$ and
$$|a|_{h^{1,\,2}_{\rm atb}(\mu)}\ls
\lf[\dint_Q|T^\ast f(x)-m_R(T^\ast f)|^{2}\,d\mu(x)\r]^{1/2}
[\mu(Q)]^{1/2}[1+\dz(Q, R)].$$
By this, \eqref{2.7}, \eqref{2.8} with $Q$ replaced by $R$
and the H\"older inequality, we see that
\begin{eqnarray*}
&&\dint_Q|T^\ast f(x)-m_R(T^\ast f)|^{2}\,d\mu(x)\\
&&\quad=\dint_Q[T^\ast f(x)-m_R(T^\ast f)]a(x)\,d\mu(x)\\
 &&\quad\le \Big|\dint_QT^\ast f(x)a(x)\,d\mu(x)\Big|
 +|m_R(T^\ast f)|\dint_Q|a(x)|\,d\mu(x)\\
 &&\quad\ls [1+\dz(Q, R)]\|f\|_{\cy^\ast}
 \lf[\dint_Q|T^\ast f(x)-m_R(T^\ast f)|^{2}\,d\mu(x)\r]^{1/2}
[\mu(Q)]^{1/2}.
\end{eqnarray*}
This in turn implies that
$$\lf[\frac1{\mu(Q)}\dint_Q|T^\ast f(x)-m_R(T^\ast f)|^{2}\,d\mu(x)\r]^{1/2}
\ls [1+\dz(Q, R)]\|f\|_{\cy^\ast},$$
which together with \eqref{2.9} and the H\"older inequality yields
\eqref{2.10}. Combining \eqref{2.8}, \eqref{2.9} and \eqref{2.10} implies the claim.

Let $h^{1,\,2}_{\rm fin}(\mu)$ be the set of all finite
linear combinations of all
$(2, 1)$-atomic blocks or $(2, 1)$-blocks $b$ as in \eqref{1.2}.
Then $h^{1,\,2}_{\rm fin}(\mu)$
is dense in $\lho$. On the other hand, $h^{1,\,2}_{\rm fin}(\mu)$
coincides with $L^2_c(\mu)$ as vector spaces. Then
by Remark \ref{r2.1} and the above claim, we
have that for all $g\in h^{1,\,2}_{\rm fin}(\mu)$ and $f\in
\cy^\ast$ with $\|f\|_{\cy^\ast}=1$,
\begin{eqnarray*}
\lf|\langle Tg, f\rangle\r|=\lf|\langle g, T^\ast f\rangle\r|
\ls \|g\|_{\lho}\|T^\ast f\|_\lbmo\ls \|g\|_{\lho}.
\end{eqnarray*}
This together with \eqref{1.2} implies that $Tg\in\cy$ and
$\|Tg\|_\cy\ls \|g\|_{\lho},$ which via a density argument then
completes the proof of Theorem \ref{t1.1} (ii). This finishes
the proof of Theorem \ref{t1.1}.
\end{pftf}
%

\newtheorem{pfts}{\bf  Proof of Theorem \ref{t1.2}.}
\renewcommand\thepfts{}

\begin{pfts}\rm
By similarity we only prove (i). As in the proof of Theorem \ref{t1.1},
we choose $\eta=2$ and $\gz=1$ in the definitions of $\hpg$.
Let $f\in H^1(\mu)$ and $f=\sum_{i=1}^\fz b_i$,
where for each $i\in\nn$, $b_i$ is a $(p, 1)$-atomic block with $p$
as in the theorem.
Since $H^1(\mu)\subset \lo$ and $T$ is bounded from $\lo$ to $\wlo$,
we see that $Tf$ is well defined. Furthermore,
by the boundedness of $T$ from $\lo$ to $\wlo$, we have that for any $\ez>0$,
\begin{eqnarray*}
\lim_{N\to\fz}\mu\lf(\lf\{x\in\rd:
\lf|T\lf(\sum_{i=N+1}^\fz b_i\r)(x)\r|>\ez\r\}\r)
\ls \lim_{N\to\fz}\frac1\ez\sum_{i=N+1}^\fz\|b_i\|_\lo=0.
\end{eqnarray*}
This via the Riesz theorem implies that
there exists a subsequence
$\{T(\sum_{i=1}^{j_k} b_i)\}_{j_k}$ of $\{T(\sum_{i=1}^j b_i)\}_j$
such that for $\mu$-a.\,e.\, $x\in\rd$,
\begin{eqnarray*}
|Tf(x)|\le \lf|T\lf(\sum_{i=1}^{j_k-1}b_i\r)(x)\r|+
\lf|T\lf(\sum_{i=j_k}^\fz b_i\r)(x)\r|
&\le&\sum_{i=1}^{j_k-1}|Tb_i(x)|+\lf|T\lf(\sum_{i=j_k}^\fz b_i\r)(x)\r|\\
&\to& \sum_{i=1}^\fz|Tb_i(x)|,\, j_k\to\fz.
\end{eqnarray*}
Since $T$ is sublinear, then from this fact,
we deduce that for $\mu$-a.\,e.\, $x\in\rd$,
$|Tf(x)|\ls \sum_{i=1}^\fz |Tb_{i}(x)|$,
which together with \eqref{1.1} in turn implies that
\begin{equation*}
\|Tf\|_\lo\ls \sum_{i=1}^\fz \|Tb_{i}\|_\lo\ls \sum_{i=1}^\fz|b_i|_\hp.
\end{equation*}
By this, we have that $Tf\in\lo$ and $\|Tf\|_\lo\ls\|f\|_{H^1(\mu)}.$
This finishes the proof of Theorem \ref{t1.2} (i), and hence,
the proof of Theorem \ref{t1.2}.
\end{pfts}

\section{Applications\label{s3}}

\hskip\parindent
In this section, we apply Theorems \ref{t1.1} to
Calder\'on-Zygmund operators, Riesz potentials and multilinear commutators
generated by Calder\'on-Zygmund operators
or fractional integral operators with Lipschitz functions,
to simplify the existing proofs in the corresponding papers.
We also prove that if $\rd$ is an initial cube,
then the Calder\'on-Zygmund operator is bounded from $\lho$ to $\lo$.

\subsection{Calder\'on-Zygmund operators and Riesz potentials\label{s3.1}}

\hskip\parindent
Recall that a $\mu$-locally integrable function $K$ on
$\rd\times\rd\setminus\{(x, y)\in\rd\times\rd:\,x=y\}$ is called
a Calder\'on-Zygmund kernel if there exists a positive
constant $C$ such that for all $x$, $y\in\rd$ with $x\not=y$,
\begin{equation}\label{3.1}
 |K(x, y)|\le C\frac 1{|x-y|^n},
\end{equation}
and for all $x$, $x'$ and $y\in\rd$ with $|x-x'|\le |x-y|/2$,
\begin{equation}\label{3.2}
|K(x, y)-K(x', y)|+|K(y, x)-K(y, x')|\le C\frac{|x-x'|}{|x-y|^{n+1}}.
\end{equation}
For all $\ez\in(0, \fz)$, $x\in\rd$ and $f\in\lt$,
define the truncated operator $T_\ez$ by
\begin{equation*}
T_\ez f(x)\equiv\int_{|x-y|\ge\ez} K(x, y)f(y)\,d\mu(y).
\end{equation*}
It is known that if the operators $\{T_\ez\}_{\ez>0}$
are bounded on $\lt$ uniformly, then there exists an operator $T$ which is the
weak limit as $\ez\to0$ of certain subsequence of
$\{T_\ez\}_{\ez>0}$; see \cite{t01}. It was proved in \cite{cmy05}
that the operator $T$ is also bounded on $\lt$ and
satisfies that for all $f\in\lt$ with bounded support
and $x\notin\supp(f)$,
\begin{equation}\label{3.3}
Tf(x)\equiv\int_\rd K(x, y)f(y)\,d\mu(y).
\end{equation}

The following Proposition \ref{p3.1} was claimed in \cite{cmy05}
without a proof. Using Theorem \ref{t1.1}, we can give
a simpler proof of Proposition \ref{p3.1} as below.

\begin{prop}\label{p3.1}
Let $T$ be a bounded linear operator on $\lt$ as in \eqref{3.3}
with the kernel $K$ satisfying \eqref{3.1} and \eqref{3.2}.
Then $T$ extends to a bounded linear operator from $H^1(\mu)$ to $\lo$.
\end{prop}

\begin{pf}\rm
Let $b\equiv \lz_1a_1+\lz_2 a_2$ be any $(2, 1)$-atomic block as
in Definition \ref{d1.4}.
Since $T$ is linear, we write
\begin{eqnarray*}
\|Tb\|_\lo&&\le\sum_{j=1}^2|\lz_j|\dint_{2Q_j}|Ta_j(x)|\,d\mu(x)
+\sum_{j=1}^2|\lz_j|\dint_{(2\sqrt{d}R)\setminus (2Q_j)}\cdots\\
&&\hs+\dint_{\rd\setminus (2\sqrt{d}R)}|Tb(x)|\,d\mu(x)
\equiv{\rm I_1}+{\rm I_2}+{\rm I_3},
\end{eqnarray*}
where for $j=1$, $2$, $Q_j$ and $R$ are as in Definition \ref{d1.4}.
Using the H\"older inequality, the boundedness of $T$ on $\lt$ and
Definition \ref{d1.4}, we have that 
${\rm I_1}\ls \sum_{j=1}^2|\lz_j|$.
By an argument similar to that used in \cite[p.\,113-114]{t01}, we obtain
${\rm I_2}+{\rm I_3}\ls \sum_{j=1}^2|\lz_j|$, which combined the estimate of
${\rm I_1}$ implies \eqref{1.1} with
$\cy=\lo$. This together with Theorem \ref{t1.1} and
$H^1(\mu)= H^{1,\,2}_{\rm atb}(\mu)$ with equivalent norms yields the boundedness
of $T$ from $H^1(\mu)$ to $\lo$, which completes the proof
of Proposition \ref{p3.1}.
\end{pf}

Let $T$ be as in Proposition \ref{p3.1}.
Recall that $T^\ast1=0$ means that for any bounded function $b$ with compact
support and $\int_\rd b(x)\,d\mu(x)=0$,
$$\int_\rd Tb(x)\,d\mu(x)=0;$$
see \cite{cmy05}. By Proposition \ref{p3.1}, this makes sense.

If $T^\ast 1=0$, using Proposition \ref{p3.1} again, we can
also complete the proof of Theorem 1 in \cite{cmy05} as follows.

\begin{prop}\label{p3.2}
Let $T$ be the same as in Proposition \ref{p3.1} and $T^\ast 1=0$.
Then $T$ extends to a bounded linear operator on $H^1(\mu)$.
\end{prop}

\begin{pf}\rm Let $\cm_\Phi$ be as in Definition \ref{d1.2} and take
$\eta=4$ in the definition of $H^{1,\,2}_{\rm atb,\,2}(\mu)$.
For any $f\in H^{1,\,2}_{\rm atb,\,2}(\mu)$, by Definition \ref{d1.4},
there exist $(2, 2)$-atomic blocks $\{b_i\}_{i=1}^\fz$ such
that $f=\sum_{i=1}^\fz b_i$ and
$\sum_{i=1}^\fz|b_i|_{H^{1,\,2}_{\rm atb,\,2}(\mu)}<\fz$. By an
argument similar to that used in the proof
of Theorem 1 in \cite{hmy05}, we have that for all $b_i$,
\begin{equation}\label{3.4}
\|\cm_\Phi(Tb_i)\|_{L^1(\mu)}\ls |b_i|_{H_{{\rm atb},\,2}^{1,\,2}(\mu)}.
\end{equation}

On the other hand, by an argument similar to the proof of Theorem 4.2 in
\cite{t01}, we obtain
\begin{equation}\label{3.5}
\lf\|\sum_{i=1}^\fz|Tb_i|\r\|_\lo\le \sum_{i=1}^\fz\|Tb_i\|_\lo
\ls \sum_{i=1}^\fz|b_i|_{H^{1,\,2}_{\rm atb,\,2}(\mu)}<\fz.
\end{equation}
Observe that for each $x\in\rd$ and $\vz\sim x$,
there exists a positive constant $M$, depending on $x$, such that for all $y\in\rd$,
$0\le \vz(y)\le M$. Moreover, by Proposition \ref{p3.1},
we obtain that $Tf=\sum_{i=1}^\fz Tb_i$ in $\lo$.
These two facts together with \eqref{3.5} and
the Lebesgue dominated convergence theorem yield that
\begin{eqnarray*}
\dint_{\rd}\vz(y)Tf(y)\,d\mu(y)=\dint_\rd\sum_{i=1}^\fz\vz(y)Tb_i(y)\,d\mu(y)
=\sum_{i=1}^\fz\dint_\rd\vz(y)Tb_i(y)\,d\mu(y).
\end{eqnarray*}
From this, it further follows that for all $x\in\rd$,
$$\cm_\Phi(Tf)(x)\le \sum_{i=1}^\fz\cm_\Phi(Tb_i)(x),$$
which together with the Levi lemma and \eqref{3.4} yields that
$$\|\cm_\Phi(Tf)\|_\lo\le \lf\|\sum_{i=1}^\fz\cm_\Phi(Tb_i)\r\|_\lo\le
\sum_{i=1}^\fz\|\cm_\Phi(Tb_i)\|_\lo\ls
\sum_{i=1}^\fz |b_i|_{H_{{\rm atb},\,2}^{1,\,2}(\mu)}.$$
This together with Definition \ref{d1.2} and $H^1(\mu)= H^{1,\,2}_{\rm atb,\,2}(\mu)$
with equivalent norms in turn implies that $Tf\in H^1(\mu)$ and
$\|Tf\|_{H^1(\mu)}\ls \|f\|_{H^1(\mu)}$,
which completes the proof of Proposition \ref{p3.2}.
\end{pf}

\begin{prop}\label{p3.3}
Let $\rd$ be an initial cube and
$T$ as in Proposition \ref{p3.1}. Then $T$ extends to a bounded linear
operator from $\lho$ to $\lo$.
\end{prop}

\begin{pf}\rm
By Theorem \ref{t1.1} (ii),
we only need to prove that $T$ satisfies \eqref{1.2} with $\cy=\lo$.
Following the proof of Proposition \ref{3.1}, we see
that for all $(2, 1)$-atomic blocks $b$ with $\supp(b)\subset R$ and $R\notin\dm$,
$\|Tb\|_\lo\ls |b|_{h^{1,\,2}_{\rm atb}(\mu)}.$
Now assume that $b$ is a $(2, 1)$-block with $\supp(b)\subset R$ and $R\in\dm$.
Since $T$ is linear, we write
$\|Tb\|_\lo\le {\rm I_1}+{\rm I_2}+{\rm I_3},$
where ${\rm I_j}$, $j=1,\,2,\,3$, are as in Proposition \ref{p3.1}.
By the same argument as in the proof of Proposition \ref{p3.1},
we have that
${\rm I_1}+{\rm I_2}\ls |b|_{h_{\rm atb}^{1,\,2}(\mu)}$.
It remains to estimate ${\rm I_3}$. Since $R\in\dm$ and $\rd$ is an initial cube,
we see that $\dz(R, \rd)\ls 1$.
From this together with  \eqref{3.1}, \eqref{3.3}, Definition \ref{d1.6}
and the fact that for all $x\in \rd\setminus (2\sqrt{d}R)$
and $y\in R$, $|x-x_R|\ls |x-y|$, it follows that
\begin{eqnarray*}
{\rm I_3}&\ls& \dint_{\rd\setminus (2\sqrt{d}R)}\dint_R\frac{|b(y)|}
{|x-y|^n}\,d\mu(y)\,d\mu(x)\\
&\ls& \dint_{\rd\setminus (2\sqrt{d}R)}\frac{\|b\|_\lo}{|x-x_R|^n}\,d\mu(x)
\ls \dz(R, \rd)\|b\|_\lo\ls |b|_{h_{\rm atb}^{1,\,2}(\mu)}.
\end{eqnarray*}
This combined with the estimates of ${\rm I_1}$ and ${\rm I_2}$
implies for all $(2, 1)$-blocks $b$ with $\supp(b)\subset R$ and $R\in\dm$,
$\|Tb\|_\lo\ls |b|_{h^{1,\,2}_{\rm atb}(\mu)},$
which together with the estimate for $(2, 1)$-atomic blocks $b$ with
$\supp(b)\subset R$ and $R\notin\dm$ further completes
the proof of Proposition \ref{p3.3}.
\end{pf}

We now consider Riesz potentials in \cite{gg04}.
Let $\az\in(0, n)$ and $K_\az$ be a locally integrable function on
$\rd\times\rd\setminus \{x=y\}$ satisfying
that there exists a positive
constant $C$ such that for all $x$, $y\in\rd$ with $x\not=y$,
\begin{equation}\label{3.6}
 |K_\az(x, y)|\le C\frac 1{|x-y|^{n-\az}},
\end{equation}
and for all $x$, $x'$ and $y\in\rd$ with $|x-x'|\le |x-y|/2$,
\begin{equation}\label{3.7}
|K_\az(x, y)-K_\az(x', y)|+|K_\az(y, x)-K_\az(y, x')|
\le C\frac{|x-x'|^\dz}{|x-y|^{n-\az+\dz}},
\end{equation}
where $\dz\in(0, 1]$. The operator
$T^\az$ associated to the above kernel $K_\az$
and the measure $\mu$ is defined by setting,
for all $f\in\lt$ with bounded support and
$x\notin \supp(f)$,
\begin{equation}\label{3.8}
T^\az f(x)\equiv\int_{\rd} K_\az(x, y)f(y)\,d\mu(y).
\end{equation}

The operator $T^\az$ was introduced by Garc\'ia-Cuerva and Gatto in \cite{gg04}.
By sealing the gap existing in the proof that (III) implies (IV) of
Theorem 1.1 in \cite{hmy08}, we have the boundedness
of $T^\az$ as follows.

\begin{prop}\label{p3.4}
Let $\az\in(0, n)$ and $T^\az$ be a linear operator as in \eqref{3.8}
with the kernel $K_\az$ satisfying \eqref{3.6} and \eqref{3.7}.
Then $T^\az$ extends to a bounded linear operator from $H^1(\mu)$
to $L^{n/(n-\az)}(\mu)$.
\end{prop}

\begin{pf}\rm
Take $\eta=4$ in the definition of $H^{1,\, n/\az}_{\rm atb}(\mu)$ and
let $b\equiv \lz_1a_1+\lz_2 a_2$ be any $(n/\az, 1)$-atomic block.
Since $T^\az$ is linear, we write
\begin{eqnarray*}
\|T^\az b\|_{L^{n/(n-\az)}(\mu)}^{n/(n-\az)}
&&\ls\sum_{j=1}^2|\lz_j|^{n/(n-\az)}\dint_{2Q_j}
|T^\az a_j(x)|^{n/(n-\az)}\,d\mu(x)
+\sum_{j=1}^2|\lz_j|^{n/(n-\az)}\\
&&\hs\times\dint_{(2\sqrt{d}R)\setminus (2Q_j)}\cdots
+\dint_{\rd\setminus (2\sqrt{d}R)}|T^\az b(x)|^{n/(n-\az)}\,d\mu(x)
\equiv{\rm L_1}+{\rm L_2}+{\rm L_3},
\end{eqnarray*}
where for $j=1$, $2$, $Q_j$ and $R$ are as in Definition \ref{d1.4}.
Recall that $T^\az$ is bounded from $\lp$ to $L^q(\mu)$ for all
$p\in (1, n/\az)$ and $q$ with $1/q=1/p-n/\az$ (see \cite{gg04}).
By an argument similar to the proof in \cite[pp.\,376-380]{hmy08},
we have that for all cubes $Q$ and functions $a\in L^{n/\az}(\mu)$ supported
in $Q$,
$$\dint_Q|T^\az a(x)|^{n/(n-\az)}\,d\mu(x)\ls
\|a\|_{L^{n/\az}(\mu)}^{n/(n-\az)}\mu(2Q).$$
From this and Definition \ref{d1.4}, it follows that
${\rm L_1}\ls \sum_{j=1}^2|\lz_j|^{n/(n-\az)}$.
Moreover, arguing as the proof in \cite[p.\,381]{hmy08}, we obtain  that
${\rm L_2}+{\rm L_3}\ls \sum_{j=1}^2|\lz_j|^{n/(n-\az)}$.
This together with the estimate of ${\rm L_1}$ implies \eqref{1.1}
with $\cy=L^{n/(n-\az)}(\mu)$,
from which, Theorem \ref{t1.1} and $H^1(\mu)= H^{1,\,n/\az}_{\rm atb}(\mu)$
with equivalent norms, it follows that
$T^\az$ extends to a bounded linear operator from $H^1(\mu)$ to $L^{n/(n-\az)}(\mu)$.
This finishes the proof of Proposition \ref{p3.4}.
\end{pf}

\subsection{Multilinear commutators}\label{s3.2}

\hskip\parindent This subsection is devoted to the boundedness of
multilinear commutators of Calder\'on-Zygmund operators
and fractional integral operators
with Lipschitz functions. We begin with the definition of
Lipschitz functions in \cite{gg05}.

\begin{defn}\label{d3.1}\rm
 Let $\bz\in(0, \fz)$. A function $f\in L^1_{\loc}(\mu)$ is said to
 belong to the space $Lip\,(\bz, \mu)$ if there exists a positive constant
 $C$ such that for $\mu$-almost every $x$ and $y\in\supp(\mu)$,
\begin{equation}\label{3.9}
|f(x)-f(y)|\le C|x-y|^\bz.
\end{equation}
Moreover, we define the $Lip\,(\bz, \mu)$ norm of $f$ to be the minimal constant
$C$ in \eqref{3.9} and denote it by $\|f\|_{Lip\,(\bz,\,\mu)}$.
\end{defn}

Let $T$ be a bounded linear operator on $\lt$ as in \eqref{3.3}
with the kernel $K$ satisfying \eqref{3.1} and \eqref{3.2},
$m\in\nn$, $\bz_i\in(0, 1]$ and $h_i\in Lip\,(\bz_i, \mu)$, $i=1,\,\cdots,\, m$.
The multilinear commutator $T_{\vec h}$ is formally defined by
\begin{equation}\label{3.10}
T_{\vec h}(f)\equiv[h_m,\,\cdots,[h_2,\,[h_1, T]]\cdots](f),
\end{equation}
where $\vec{h}\equiv(h_1,\,h_2,\cdots,\,h_m)$ and
\begin{equation*}
[h_1, T](f)\equiv h_1T(f)-T(h_1f).
\end{equation*}

The operator $T_{\vec h}$ was introduced in \cite{my06}
and the following Proposition \ref{p3.5} was
also obtained there (see \cite[Theorem 3.1]{my06}).
Using Theorem \ref{t1.1},
we can also give a simpler proof of this proposition
as below.

\begin{prop}\label{p3.5}
Let $m\in\nn$, $\bz_i\in(0, 1]$,
$h_i\in Lip\,(\bz_i, \mu)$ for $i=1, \,\cdots,\,m$,
and $T_{\vec h}$ be as in \eqref{3.10}.
If $\bz\equiv\sum_{i=1}^m\bz_i<n$ and
$1/q=1-\bz/n$,
then $T_{\vec h}$ extends to a bounded linear operator from
$H^1(\mu)$ to $L^q(\mu)$.
\end{prop}

\begin{pf}\rm
Take $\eta =4$
in the definition of $H^{1,\, n/\bz}_{\rm atb}(\mu)$.
Repeating the proof of Theorem 3.1 in
\cite{my06}, we have that for all
$(n/\bz, 1)$-atomic blocks $b$,
$\|T_{\vec h}(b)\|_{L^q(\mu)}\ls |b|_{H^{1,\,n/\bz}_{\rm atb}(\mu)}$.
This implies \eqref{1.1} with $\cy=L^q(\mu)$. Since $T_{\vec h}$ is linear,
then an application of Theorem \ref{t1.1} together with
$H^1(\mu)= H^{1,\,n/\bz}_{\rm atb}(\mu)$
with equivalent norms yields that
$T_{\vec h}$ extends to a bounded linear operator from
$H^1(\mu)$ to $L^q(\mu)$, which
completes the proof of Proposition \ref{p3.5}.
\end{pf}

We now consider multilinear commutators
generated by fractional integral operators and
Lipschitz functions in \cite{my06}. To be precise,
let $\az\in(0, n)$, $x\in\supp(\mu)$ and $f\in\lin$ with bounded support.
The fractional integral operator $I_\az$ is defined by
$$I_\az(f)(x)\equiv\dint_\rd\frac{f(y)}{|x-y|^{n-\az}}\,d\mu(y).$$
In \cite{gm01}, Garc\'ia-Cuerva and Martell introduced the operator $I_\az$
and proved that $I_\az$ is bounded from $\lp$ to $L^{q,\,\fz}(\mu)$ with
$p\in[1, n/\az)$ and $1/q=1/p-\az/n$.
For $m\in\nn$ and $h_i\in Lip\,(\bz_i, \mu)$,
where $\bz_i\in(0, 1]$, $i=1,\,\cdots,\,m$, and $\az+\sum_{i=1}^m\bz_i<n$, define
the multilinear commutator
$I_{\az,\,\vec h}$ by setting, for all $f\in\lt$ with bounded support
and $x\notin \supp(f)$,
\begin{equation}\label{3.11}
I_{\az,\,\vec h}(f)(x)\equiv\dint_\rd\prod_{i=1}^m[h_i(x)-h_i(y)]
\frac{f(y)}{|x-y|^{n-\az}}\,d\mu(y).
\end{equation}

The commutator $I_{\az,\,\vec h}$ was also introduced in \cite{my06}
and the following Proposition \ref{p3.6} is Theorem 4.2 in \cite{my06}.
Applying Theorem \ref{t1.1}, we can also give a simpler proof
of Proposition \ref{p3.6} as below.

\begin{prop}\label{p3.6}
Let  $\az\in(0, n)$, $m\in\nn$,
$\bz_i\in(0, 1]$, $h_i\in Lip\,(\bz_i, \mu)$
for $i=1, \,\cdots,\,m$,
and $I_{\az,\,\vec h}$ be as in \eqref{3.11}.
If $\bz\equiv\az+\sum_{i=1}^m\bz_i<n$ and
$1/q=1-\bz/n$,
then $I_{\az,\,\vec h}$ extends to a bounded linear operator from
$H^1(\mu)$ to $L^q(\mu)$.
\end{prop}

\begin{pf}\rm
Take $\eta =4$ in the definition of $H^{1,\, n/\bz}_{\rm atb}(\mu)$.
Arguing as in the proof of Theorem 3.1 in \cite{my06}, we see that for all
$(n/\bz, 1)$-atomic blocks $b$,
$\|I_{\az,\,\vec h}(b)\|_{L^q(\mu)}\ls |b|_{H^{1,\,n/\bz}_{\rm atb}(\mu)}$,
which implies \eqref{1.1} with $\cy=L^q(\mu)$.
From this together with Theorem \ref{t1.1} and $H^1(\mu)
= H^{1,\,n/\bz}_{\rm atb}(\mu)$
with equivalent norms, it follows that
$I_{\az,\,\vec h}$ extends to a bounded linear operator
from $H^1(\mu)$ to $L^q(\mu)$. This
finishes the proof of Proposition \ref{p3.6}.
\end{pf}

\bigskip

\noindent Dachun Yang and  Dongyong Yang

\medskip

\noindent School of Mathematical Sciences,
 Beijing Normal University, Laboratory of Mathematics and Complex systems,
Ministry of Education, Beijing 100875, People's Republic of China

\medskip

\noindent{\it E-mail addresses}: \texttt{dcyang@bnu.edu.cn}

\hspace{2.56cm}\texttt{dyyang@mail.bnu.edu.cn}
\end{document}